\documentclass[final,1p,times]{elsarticle}

\usepackage{amssymb}
 \usepackage{amsthm}
\usepackage{amscd}
\usepackage{amsmath}
\usepackage{amsfonts}
\usepackage{amssymb}
\usepackage{graphicx}
\newtheorem{theorem}{Theorem}

\newtheorem{lemma}[theorem]{Lemma}

\usepackage{mathrsfs}
\usepackage{titletoc}

\newcommand{\ra}{\rightarrow}
\newcommand{\p}{\partial}
\newcommand{\f}{\frac}

\newcommand{\be}{\begin{equation}}
\renewcommand{\ra}{\rightarrow}
\newcommand{\ee}{\end{equation}}
\newcommand{\bea}{\begin{eqnarray}}
\newcommand{\eea}{\end{eqnarray}}
\newcommand{\bna}{\begin{eqnarray*}}
\newcommand{\ena}{\end{eqnarray*}}

\renewcommand{\le}{\left}
\newcommand{\ri}{\right}

\journal{***}

\begin{document}

\begin{frontmatter}
\title{Calculus of variations on locally finite graphs}

\author{Yong Lin}
\ead{yonglin@mail.tsinghua.edu.cn}
\address{Yau Mathematical Sciences Center, Tsinghua University, Beijing 100084, China}

\author{Yunyan Yang}
\ead{yunyanyang@ruc.edu.cn}
\address{ Department of Mathematics,
Renmin University of China, Beijing 100872, China}

\begin{abstract}
Let $G=(V,E)$ be a locally finite graph. Firstly, using calculus of variations, including a direct method of
variation and the mountain-pass theory, we get sequences of solutions to several local equations on $G$ (the Schr\"odinger equation,
the mean field equation, and the Yamabe equation). Secondly, we derive uniform estimates for those local solution sequences.
Finally, we obtain global solutions by extracting convergent sequence of solutions. Our method can be
described as a variational method from local to global.
\end{abstract}


\end{frontmatter}

\section{Introduction}
Partial differential equations on Euclidean space or manifolds are important topics in mathematical physics and differential geometry.
As their discrete versions, it is important to study the difference equations on graph, particularly the existence problem for
such equations.

About five years ago, joined with Grigor'yan, we
 systematically raised and studied Kazdan-Warner equations, Yamabe equations and Sch\"ordinger equations on graphs in \cite{GLY1,GLY2,GLY3}.
 We first established the Sobolev spaces and the functional framework. Then the problem of solving the equations is transformed into
 finding critical points of various functionals. As a consequence, variational methods are applied to these problems. If the graph is finite,
 then all the Sobolev spaces have finite dimensions, and whence they are pre-compact. For this reason, the variational problems
 for finite graph are comparatively simple \cite{GLY1,GLY2}. Since the graph has no concept of dimension, if it includes infinite
 vertices,
 the Sobolev embedding theorems becomes unusual. An easy-to-understand one was observed by us \cite{GLY3} under the assumption that
 the graph is locally finite and its measure has positive lower bound (see next section for details).
 Any other Sobolev embedding theorem for infinite graph would be extremely interesting.

 In recent years, the research in this field has aroused great interest. Motivated by \cite{YangJFA12,GLY3},
 Zhang-Zhao \cite{Zhang-Zhao} obtained  nontrivial solutions to certain nonlinear Schr\"odinger equation. Similar equations on infinite metric graphs were studied by Akduman-Pankov \cite{Akduman-Pankov}.
 The Kazdan-Warner equation was extended by Keller-Schwarz \cite{Keller-Schwarz} to canonically compactifiable graphs, and by Ge-Jiang
 \cite{Ge-Jiang} to certain infinite graph. For other related works, we refer the reader to \cite{Huang-Lin-Yau,Liu-Yang,Man,Han-Shao-Zhao,
 Tian-Zhang-Zhang,Zhang-Lin} and the references therein.

 In this paper, we study various equations on locally finite graphs, say Schr\"ordinger equation, Mean field equation and Yamabe equation.
 Assuming that the weights of the graph have a positive lower bound and the distance function of the graph belongs to $L^p$, we derive a Sobolev embedding theorem, which is crucial in our analysis. In addition to
 the Sobolev embedding theorem, we also employ calculus of variations, including a direct method of variation and the mountain-pass theorem.
 It is remarkable that we show how to get solutions from local to global by using variational method.

\section{Notations and main results}\label{sec1}

Let $G=(V,E)$ be a connected graph, where $V$ denotes the vertex set and $E$ denotes the edge set.
For any edge $xy\in E$, we assume that its weight $w_{xy}>0$ and that $w_{xy}=w_{yx}$.
Let $\mu:V\ra \mathbb{R}^+$ be a finite measure. For any function $u:V\ra \mathbb{R}$, the Laplacian
of $u$ is defined as
\be\label{lap}\Delta u(x)=\f{1}{\mu(x)}\sum_{y\sim x}w_{xy}(u(y)-u(x)),\ee
where $y\sim x$ means $xy\in E$ or $y$ is adjacent to $x$.
The gradient form is written by
\be\label{grad-form}\Gamma(u,v)(x)=\f{1}{2\mu(x)}\sum_{y\sim x}w_{xy}(u(y)-u(x))(v(y)-v(x)).\ee
Denote $\Gamma(u)=\Gamma(u,u)$ and $\nabla u\nabla v=\Gamma(u,v)$. The length of the gradient of $u$ is represented by
\be\label{grd}|\nabla u|(x)=\sqrt{\Gamma(u)(x)}=\le(\f{1}{2\mu(x)}\sum_{y\sim x}w_{xy}(u(y)-u(x))^2\ri)^{1/2}.\ee
The integral of a function $f$ on $V$ is given as
\be\label{int}\int_V fd\mu=\sum_{x\in V}\mu(x)f(x).\ee
For any $q>0$, we let $L^q(V)$ be a linear space of functions $f:V\ra\mathbb{R}$ with the norm
\be\label{Lp}\|f\|_{L^q(V)}=\le(\int_V|f|^qd\mu\ri)^{1/q}.\ee
While $L^\infty(V)$ includes all functions $f:V\ra\mathbb{R}$ satisfying
$$\|f\|_{L^\infty(V)}=\sup_{x\in V}|f(x)|<\infty.$$

If $x,y\in V$ and $y$ is adjacent to $x$, then the distance between $x$ and $y$ is defined as $1$. While if $y$ is not adjacent to $x$,
then there exists a shortest path $\gamma$ connecting $y$ and $x$, and thus the distance between $x$ and $y$ is defined as the number
of edges belonging to $\gamma$. Given any $O\in V$. Denote
the distance between $x$ and $O$ by
\be\label{distance}\rho(x)=\rho(x,O).\ee  For any integer $k\geq 1$, we denote a ball centered at $O$ with radius $k$ by
\be\label{Ball-k}B_k=B_k(O)=\le\{x\in V: \rho(x)<k\ri\}.\ee
The boundary of $B_k$ is written as
\be\label{boundary}\p B_k=\{x\in V:\rho(x)=k\}.\ee
According to \cite{GLY2}, $W_0^{1,2}(B_k)$ stands for a Sobolev space including all functions $u:B_k\ra\mathbb{R}$ with $u=0$
on the boundary $\p B_k$ given as in (\ref{boundary}). For any fixed $k$, it is pre-compact. Precisely, if $(u_j)$ is a bounded sequence in
$W_0^{1,2}(B_k)$, i.e.
\be\label{W012}\|u_k\|_{W_0^{1,2}(B_k)}=\le(\int_{B_k}|\nabla u_k|^2d\mu\ri)^{1/2}\leq C,\ee
where the notations in (\ref{grad-form}), (\ref{grd}) and (\ref{int}) are used,
 then there exists a subsequence of $(u_j)$ converging to some function $u$ under the norm in (\ref{W012}).

Recall another important Sobolev space $W^{1,2}(V)$ including all functions $u:V\ra\mathbb{R}$ with
\be\label{W12}\|u\|_{W^{1,2}(V)}=\le(\int_V(|\nabla u|^2+u^2)d\mu\ri)^{1/2}<+\infty.\ee
Let $C_c(V)$ be a set of all functions with finite support, and $W_0^{1,2}(V)$ be a completion of $C_c(V)$ under the norm as in
(\ref{W12}). Both of $W^{1,2}(V)$ and $W_0^{1,2}(V)$ are Hilbert spaces with the same inner product
$\langle u,v\rangle=\int_V(\nabla u\nabla v+uv)d\mu$.

A connected graph is said to be locally finite if for any fixed $O\in V$, $B_k$ is a finite subgraph.
In \cite{GLY3}, we made a key observation under the assumption that $G$ is locally finite, and there exists a constant $\mu_0>0$ satisfying
\be\label{mu0}\mu(x)\geq \mu_0\quad{\rm for\,\,all}\quad x\in V.\ee
Namely, a Sobolev embedding theorem holds.
\begin{theorem}[\cite{GLY3}]\label{Soblev-1}
Let $G=(V,E)$ be a connected locally finite graph. If (\ref{mu0}) is satisfied,
then for any $u\in W^{1,2}(V)$ and any $2\leq q\leq\infty$, there exists a positive constant $C$ depending only on $q$ and $\mu_0$
satisfying
$\|u\|_{L^q(V)}\leq C\|u\|_{W^{1,2}(V)}$.
In particular,
$$\|u\|_{L^\infty(V)}\leq \f{1}{\sqrt{\mu_0}}\|u\|_{W^{1,2}(V)}.$$
\end{theorem}

If instead of (\ref{mu0}), there exists some constant $w_0>0$ such that
\be\label{w-lower}w_{xy}\geq w_0\quad {\rm for\,\, all}\quad y\sim x,\,\,x,y\in V, \ee
and the distance function $\rho(x)$ defined as in (\ref{distance}) belongs to $L^p(V)$, we shall prove a Sobolev embedding as follows.
\begin{theorem}\label{Soblev-2}
Let $G=(V,E)$ be a connected locally finite graph. If the weights $w_{xy}$ satisfy (\ref{w-lower}), and the distance function $\rho(x)=\rho(x,O)\in L^p(V)$ for some
$p>0$ and some $O\in V$, then there exists some constant $C$ depending only on $w_0$, $\mu(O)$ and $p$ such that
$$\|u\|_{L^p(V)}\leq C(\|\rho\|_{L^p(V)}+1)\|u\|_{W^{1,2}(V)}.$$
\end{theorem}

If a function $h:V\ra\mathbb{R}$ has a positive lower bound on $V$, then we define a subspace of $W_0^{1,2}(V)$, which is
also a Hilbert space, namely
\be\label{H}\mathscr{H}=\le\{u\in W_0^{1,2}(V): \int_V(|\nabla u|^2+hu^2)d\mu<\infty\ri\}\ee
with an inner product
\be\label{inner-p}\langle u,v\rangle_{\mathscr{H}}=\int_V(\nabla u\nabla v+huv)d\mu.\ee

The first equation we concern is the following linear Schr\"odinger equation
\be\label{Schord}\le\{\begin{array}{lll}
-\Delta u+hu=f\quad{\rm in}\quad V\\[1.5ex]
u\in\mathscr{H},\end{array}\ri.\ee
where $\Delta$ is the Laplacian operator defined as in (\ref{lap}), and $\mathscr{H}$ is defined as in (\ref{H}).
We now state the following existence result.
\begin{theorem}\label{thm3}
Let $G=(V,E)$ be a connected locally finite graph. Assume there is some constant $a_0>0$ such that $h(x)\geq a_0$ for all $x\in V$.
If one of the following three hypotheses is satisfied:
\item{
$(i)$ $f\in L^2(V)$;}
\item{$(ii)$  $\mu(x)\geq \mu_0>0$ for all $x\in V$, $f\in L^1(V)$;}
\item{$(iii)$ the weights of the graph satisfies (\ref{w-lower}), the distance function
$\rho(x)=\rho(x,O)\in L^p(V)$ for some $p>1$, $O\in V$, and  $f\in L^{p/(p-1)}(V)$,}\\[1.1ex]
 then the equation (\ref{Schord}) has a unique solution. If in addition $f\geq 0$ and $f\not\equiv 0$ on $V$,
 then $u(x)>0$ for all $x\in V$.
\end{theorem}

The second equation we concern is the mean field equation, which is also known as the Kazdan-Warnar equation,
namely
\be\label{KZ-equ}\Delta u=f-ge^u\quad{\rm in}\quad V.\ee

 \begin{theorem}\label{thm4}
 Let $G=(V,E)$ be a locally finite graph. Suppose that $g\leq f<0$ on $V$ and $g\in L^1(V)$. Then
 the equation (\ref{KZ-equ}) has a solution.
 \end{theorem}

 We remak that using a method of the heat equation, Ge-Jiang \cite{Ge-Jiang} obtained similar result as that of Theorem \ref{thm4}
 under different assumptions on $f$ and $g$. In the case $g>0$, it is not likely to find a
 nontrivial solution as in Theorem \ref{thm4} in general. The main difficulty is that $\int_V|\nabla u|^2d\mu$
 does not control $\|u\|_{W^{1,2}(V)}$ if $V$ is an infinite graph.
 However,
 it is natural to consider the following mean field equation
 \be\label{meanfield}\le\{\begin{array}{lll}
 -\Delta u+hu=\f{ge^u}{\int_Vge^ud\mu}-f\quad{\rm in}\quad V\\[1.2ex]
 u\in\mathscr{H}\cap L^\infty(V),
 \end{array}
 \ri.\ee
 where $h$ has a positive lower bound, and $\mathscr{H}$ is defined as in (\ref{H}).
 To seek solutions of (\ref{meanfield}), we need certain Trudinger-Moser embedding. It suffices to assume (\ref{mu0}) for
 the graph in order to get that kind of embedding. Precisely we have the following:
 \begin{theorem}\label{thm5}
 Let $G=(V,E)$ be a connected locally finite graph. Suppose  (\ref{mu0}) is satisfied, there exists some constant
 $a_0>0$ such that $h(x)\geq a_0$ for all $x\in V$, $g\geq 0$ and $g\not\equiv 0$ on $V$, $g\in L^1(V)$, and $f\in L^q(V)$ for some $q\in[1,2]$. Then
 the equation (\ref{meanfield})
 has a solution.
 \end{theorem}
 Note that in Theorem \ref{thm5}, the function $f$ allows the form $\sum_{i=1}^\ell c_i\delta_{x_i}$ for some
 constants $c_1,\cdots,c_\ell$, where $\delta_{x_i}$ stands for the Dirac function satisfying
 $$\delta_{x_i}(x)=\le\{\begin{array}{lll}
 1&{\rm if}& x=x_i\\[1.5ex]
 0&{\rm if}& x\not=x_i.
 \end{array}
 \ri.$$
 As a consequence, it makes sense to consider Chern-Simons-Higgs model in locally finite graph.
 Such a model in finite graph was recently studied by Huang-Lin-Yau \cite{Huang-Lin-Yau}.

 The third equation we are interested in is the Yamabe equation
 \be\label{Yamabe}\le\{\begin{array}{lll}-\Delta u+hu=|u|^{q-2}u\quad{\rm in}\quad V\\[1.5ex]
 u\in \mathscr{H},\end{array}\ri.\ee
 where $h$ has a positive lower bound, $\mathscr{H}$ is defined as in (\ref{H}), and $q>2$. In order to find a solution
 to the equation (\ref{Yamabe}), we seek the Sobolev embedding theorem, say Theorem \ref{Soblev-1} or Theorem \ref{Soblev-2}.
 Inspired by \cite{doo,Adi-Yang,YangJFA12}, we have solved this problem in \cite{GLY3} by employing Theorem \ref{Soblev-1}.
 For application of Theorem \ref{Soblev-2}, we state the following:

 \begin{theorem}\label{thm6}
 Let $G=(V,E)$ be a connected locally finite graph.
 Let $O$ be a fixed point of $V$, the distance function
 $\rho(x)=\rho(x,O)\in L^p(V)$ for some $p>2$. Suppose  $h(x)\geq a_0>0$ for some constant $a_0$ and all $x\in V$. If further $1/h\in L^1(V)$
 or $h(x)\ra+\infty$ as $\rho(x)\ra+\infty$,
 then for any $q$ with $2<q<p$, the equation (\ref{Yamabe}) has a nontrivial solution.
 \end{theorem}

 The remaining part of this paper is organized as follows: In Section \ref{sec3}, a Sobolev embedding theorem
 (Theorem \ref{Soblev-2}) is proved; In Section \ref{sec4}, we study the linear Schr\"odinger equation, and prove Theorem \ref{thm3};
 In Section \ref{sec5}, the mean field equations are discussed, and Theorems \ref{thm4} and \ref{thm5} are proved;
 In Section \ref{sec6}, we consider the Yamabe equation and prove Theorem \ref{thm6}. Throughout this paper, we do not distinguish
 sequence and subsequence, and denote
 various constants by the same $C$.

 \section{A Sobolev embedding}\label{sec3}

 In this section, using definitions of $W^{1,2}(V)$ and $L^p(V)$, we prove Theorem \ref{Soblev-2}. \\

 {\it Proof of Theorem \ref{Soblev-2}}. Let $O$ be a fixed point in $V$. For any $x\in V$, we denote the distance
 between $x$ and $O$ by $\rho(x)=\rho(x,O)$. Choose a shortest path $\gamma=\{x_1,\cdots,x_{k+1}\}$ connecting $x$ and $O$.
 In particular $x_1=x$, $\cdots$, $x_{k+1}=O$, $x_i$ is adjacent to $x_{i+1}$ for all $1\leq i\leq k$, and $k=\rho(x)$.
 For any $u\in W^{1,2}(V)$, we get
 \be\label{tri}|u(x)|\leq |u(x_1)-u(x_2)|+\cdots+|u(x_k)-u(x_{k+1})|+|u(O)|.\ee
 Noting that (\ref{W12}) implies
 \be\label{nor}\|u\|_{W^{1,2}(V)}=\le(\sum_{z\in V,\,y\sim z}w_{zy}(u(y)-u(z))^2+\sum_{z\in V}\mu(z)u^2(z)\ri)^{1/2},\ee
 and that $\mu(z)>0$ for all $z\in V$, we have
 \be\label{uO}|u(O)|\leq \f{1}{\sqrt{\mu(O)}}\|u\|_{W^{1,2}(V)};\ee
 since $w_{zy}\geq w_0>0$ for all $z$ adjacent to $y$, in view of (\ref{nor}),
 \bea\nonumber
 \sum_{i=1}^k|u(x_i)-u(x_{i+1})|&\leq&k\max_{1\leq i\leq k}|u(x_i)-u(x_{i+1})|\\\nonumber
 &\leq&\f{k}{\sqrt{w_0}}\max_{1\leq i\leq k}\sqrt{w_{x_ix_{i+1}}}|u(x_i)-u(x_{i+1})|\\
 &\leq&\f{1}{\sqrt{w_0}}\rho(x)\|u\|_{W^{1,2}(V)}.\label{u-abs}
 \eea
 Combining (\ref{tri}), (\ref{uO}) and (\ref{u-abs}), we obtain
 \be\label{ineq}|u(x)|\leq \le(\f{1}{\sqrt{w_0}}\rho(x)+\f{1}{\sqrt{\mu(O)}}\ri)\|u\|_{W^{1,2}(V)}.\ee
  Since $\rho\in L^p(V)$ for some $p>0$ and $\rho(x,y)\geq 1$ for all $x\not =y$, in view of (\ref{Lp}), there holds
  \bna\| 1\|_{L^p(V)}=\le(\sum_{z\in V}\mu(z)\ri)^{1/p}&\leq& \le(\sum_{z\in V}\mu(z)\rho^p(z)+\mu(O)\ri)^{1/p}\\
  &\leq& 2^{1/p}\max\le\{\le(\sum_{z\in V}\mu(z)\rho^p(z)\ri)^{1/p},\,\mu(O)^{1/p}\ri\}\\
  &=&2^{1/p}\max\le\{\|\rho\|_{L^p(V)},\,\mu(O)^{1/p}\ri\}.\ena
  This together with (\ref{ineq}) leads to
  \bna
  \|u\|_{L^p(V)}&\leq&C\le(\f{1}{\sqrt{w_0}}\|\rho\|_{L^p(V)}+\f{1}{\sqrt{\mu(O)}}\|1\|_{L^p(V)}\ri)\|u\|_{W^{12}(V)}\\
  &\leq&C(\|\rho\|_{L^p(V)}+1)\|u\|_{W^{1,2}(V)}
  \ena
   for some constant $C$ depending only on $w_0$, $\mu(O)$ and $p$, as we expected. $\hfill\Box$

\section{Schr\"odinger equation}\label{sec4}

In this section, we prove Theorem \ref{thm3} by using a direct method of variation from local to global. \\

{\it Proof of Theorem \ref{thm3}}.
Fix some point $O\in V$. Denote the distance between $x$ and $O$ by $\rho(x)=\rho(x,O)$.
For any positive integer $k$, we write
$B_k=\{x\in V: \rho(x)<k\}$. Note that $h(x)\geq a_0>0$ for all $x\in V$.
Let $W_0^{1,2}(B_k)$ be the Sobolev space including all functions $u:B_k\ra \mathbb{R}$,
$u=0$ on $\p B_k$, with the norm
\be\label{norm-k}\|u\|_{W_0^{1,2}(B_k)}=\le(\int_{B_k}(|\nabla u|^2+hu^2)d\mu\ri)^{1/2}.\ee
For any fixed $k$, the norm in (\ref{norm-k}) is equivalent to that in (\ref{W012}), due to
the Poincar\'e inequality
$$\int_{B_k}u^2d\mu\leq C_k\int_{B_k}|\nabla u|^2d\mu,\quad\forall u\in W_0^{1,2}(B_k),$$
where $C_k$ is a constant depending on $k$. In general, $C_k$ tends to infinity as $k\ra\infty$.
It is convenient for us to use (\ref{norm-k}) as the norm in $W_0^{1,2}(B_k)$.
Define a functional $J_k: W_0^{1,2}(B_k)\ra \mathbb{R}$ by
\be\label{Jk}J_k(u)=\f{1}{2}\int_{B_k}(|\nabla u|^2+hu^2)d\mu-\int_{B_k}fud\mu.\ee
Set $\Lambda_k=\inf_{u\in W_0^{1,2}(B_k)}J_k(u)$. Obviously
\be\label{upbd}\Lambda_k\leq J_k(0)=0.\ee

{\it Case $(i)$. $f\in L^2(V)$.}

By the H\"older inequality and the Young inequality, we have
\bea\nonumber
\int_{B_k}fud\mu&\leq& \f{1}{\sqrt{a_0}}\le(\int_Vf^2d\mu\ri)^{1/2}\le(\int_{B_k}hu^2d\mu\ri)^{1/2}\\\nonumber
&\leq& \f{1}{\sqrt{a_0}}\|f\|_{L^2(V)}\|u\|_{W_0^{1,2}(B_k)}\\\label{f-bd}
&\leq& \f{1}{4}\|u\|_{W_0^{1,2}(B_k)}^2+\f{1}{a_0}\|f\|_{L^2(V)}^2,
\eea
where $\|u\|_{W_0^{1,2}(B_k)}$ is defined as in (\ref{norm-k}).
It follows from (\ref{Jk}) and (\ref{f-bd}) that
\be\label{J-11}J_k(u)\geq \f{1}{4}\|u\|_{W_0^{1,2}(B_k)}^2-\f{1}{a_0}\|f\|_{L^2(V)}^2.\ee
Hence
\be\label{J-12}\Lambda_k=\inf_{u\in W_0^{1,2}(B_k)}J_k(u)\geq -\f{1}{a_0}\|f\|_{L^2(V)}^2.\ee
Combining (\ref{upbd}) and (\ref{J-12}), we know that $(\Lambda_k)$ is a bounded sequence of numbers.
Now we fix a positive integer $k$ and take a sequence of functions
$(\widetilde{u}_j)\subset W_0^{1,2}(B_k)$ satisfying
\be\label{minimizing-1}J_k(\widetilde{u}_j)\ra \Lambda_k\quad{\rm as}\quad j\ra\infty.\ee
It follows from (\ref{J-11}) that $(\widetilde{u}_j)$ is bounded in $W_0^{1,2}(B_k)$. By the Sobolev embedding theorem
for finite graph \cite{GLY2}, there exists a $u_k\in W_0^{1,2}(B_k)$ such that up to a subsequence, $\widetilde{u}_j$ converges to $u_k$
under the norm (\ref{norm-k}). Clearly $J_k(u_k)=\Lambda_k$, and $u_k$ satisfies the Euler-Lagrange equation
\be\label{equ}\le\{\begin{array}{lll}
-\Delta u_k+hu_k=f&{\rm in}& B_k\\[1.5ex]
u_k=0&{\rm on}& \p B_k.
\end{array}\ri.\ee
Noting that $(\Lambda_k)$ is bounded, in view of (\ref{J-11}) and (\ref{minimizing-1}), we obtain
\be\label{bdd}\|u_k\|_{W_0^{1,2}(B_k)}^2=\int_{B_k}(|\nabla u_k|^2+hu_k^2)d\mu\leq C\ee
for some constant $C$ independent of $k$.  
For any finite set $K\subset V$, there holds $K\subset B_k$ for sufficiently large $k$.
The power of (\ref{bdd}) is evident. It ensures that
$$\|u_k\|_{L^\infty(K)}\leq \f{1}{\sqrt{a_0}\min_{x\in K}\mu(x)}\|u_k\|^2_{W_0^{1,2}(B_k)}\leq C.$$ 
Note that $(u_k)$ is naturally viewed as
a sequence of functions defined on $V$, say $u_k\equiv 0$ on $V\setminus B_k$. 
There would exist a subsequence of
$(u_k)$ (which is still denoted by $(u_k)$) and a
function $u^\ast$ such that $(u_k)$ converges to $u^\ast$ locally uniformly in $V$, i.e. for any fixed positive integer
$\ell$,
$$\lim_{k\ra\infty}u_k(x)=u^\ast(x)\quad{\rm for\,\,all}\quad x\in B_\ell.$$

Now we show that
\be\label{uast}u^\ast\in \mathscr{H}.\ee
 Since $u_k$ is viewed as a function on the whole $V$, $u_k=0$ on $V\setminus B_k$,
 and the weights of the graph is symmetric, i.e. $w_{xy}=w_{yx}$ for all $y$ adjacent to $x$,
 we have the following estimate
 \bea\nonumber
 \|u_k\|_{\mathscr{H}}^2&=&\sum_{y\sim x}w_{xy}(u_k(y)-u_k(x))^2+\sum_{x\in V}\mu(x)h(x)u_k^2(x)\\\nonumber
 &=& \sum_{y\sim x,\,x\in B_k}w_{xy}(u_k(y)-u_k(x))^2+\sum_{x\in B_k}\mu(x)h(x)u_k^2(x)\\\nonumber
 &&+\sum_{y\sim x,\,x\in \p B_k}w_{xy}(u_k(y)-u_k(x))^2\\\nonumber
 &\leq&2\sum_{y\sim x,\,x\in B_k}w_{xy}(u_k(y)-u_k(x))^2+\sum_{x\in B_k}\mu(x)h(x)u_k^2(x)\\\label{H-norm}
 &\leq&2\|u_k\|^2_{W_0^{1,2}(B_k)}.
 \eea
 Up to a subsequence, we assume $(u_k)$ converges to $u^\ast$ locally uniformly in $V$. In view of (\ref{bdd}) and (\ref{H-norm}),
 we know that $(u_k)$ is bounded in $\mathscr{H}$. Since every Hilbert space is weakly compact,
  it follows that up to a subsequence, $(u_k)$ converges to some function
 $u_1^\ast$ weakly in $\mathscr{H}$. This in particular implies 
 $$\int_Vu_k\phi d\mu\ra \int_V u_1^\ast \phi d\mu,\quad\forall \phi\in C_c(V).$$
 Let $z\in V$ be any fixed point. In the above estimate, we take $\phi$ satisfying $\phi(x)=1$ at $x=z$ and $\phi(x)=0$ at $x\not=z$. 
 Then $u_k(z)\ra u_1^\ast(z)$. Hence by the uniqueness of the limit, $u_1^\ast(z)\equiv u^\ast(z)$ for all $z\in V$, and (\ref{uast}) follows immediately.

It then follows from (\ref{equ}) that for any fixed $x\in V$, there holds
$$-\Delta u^\ast(x)+h(x)u^\ast(x)=f(x).$$
Therefore $u^\ast$ is a solution of (\ref{Schord}). To prove that $u^\ast$ is a unique solution of (\ref{Schord}),
it suffices to show the homogenueous equation
\be\label{homo}\le\{\begin{array}{lll}-\Delta u+hu=0\\[1.5ex]
u\in\mathscr{H}\end{array}\ri.\ee
has only one solution $u\equiv 0$. Since $u\in\mathscr{H}$, there exists a sequence $(\varphi_k)\subset C_c(V)$ such that
$\varphi_k\ra u$ in $\mathscr{H}$. Testing (\ref{homo}) by $\varphi_k$, we have by integration by parts
$$\langle \varphi_k,u\rangle_{\mathscr{H}}=\int_V(\nabla u\nabla\varphi_k+hu\varphi_k)d\mu=0,$$
where $\langle\cdot,\cdot\rangle_{\mathscr{H}}$ is the inner product in $\mathscr{H}$ defined as in (\ref{inner-p}).
Passing to the limit $k\ra\infty$, we conclude $\langle u,u\rangle_{\mathscr{H}}=0$, and thus $u\equiv 0$.
This confirms the uniqueness of $u^\ast$.

If $f(x)\geq 0$ for all $x\in V$, then applying the maximum principle to (\ref{equ}), we obtain $u_k(x)\geq 0$ for all
$x\in B_k$. Indeed, suppose there exists some $x_0\in B_k$ satisfying $\min_{B_k}u_k=u_k(x_0)<0$, we have by (\ref{equ}) that
$$-\Delta u_k(x_0)=f(x_0)-h(x_0)u_k(x_0)>0.$$
This is impossible, and leads to $u_k\geq 0$ on $B_k$. As a consequence, $u^\ast(x)\geq 0$ for all $x\in V$.
Since $f\not\equiv 0$, one has $u^\ast\not\equiv 0$. We now prove
$u^\ast(x)>0$ for all $x\in V$. Suppose not, there would be a point $x^\ast\in V$ such that $u^\ast(x^\ast)=0=\min_Vu^\ast$
and $\Delta u^\ast(x^\ast)>0$. It follows that
$$0>-\Delta u^\ast(x^\ast)=f(x^\ast)\geq 0,$$
which is a contradiction, and implies $u^\ast(x)>0$ for all $x\in V$. \\

{\it Case $(ii)$. $\mu(x)\geq \mu_0>0$ for all $x\in V$.}

By the Sobolev embedding theorem (Theorem \ref{Soblev-1}), we have for all $u\in W_0^{1,2}(B_k)$,
$$\|u\|_{L^\infty(B_k)}\leq \f{1}{\sqrt{\mu_0}}\|u\|_{W_0^{1,2}(B_k)}.$$
Similar to (\ref{f-bd}), there holds
\bna\int_{B_k}fud\mu&\leq& \|u\|_{L^\infty(B_k)}\|f\|_{L^1(B_k)}\\
&\leq& \f{1}{\sqrt{\mu_0}}\|u\|_{W_0^{1,2}(B_k)}\|f\|_{L^1(V)}\\
&\leq&\f{1}{4}\|u\|_{W_0^{1,2}(B_k)}^2+\f{1}{\mu_0}\|f\|_{L^1(V)}^2.\ena
In the same way,  for any $u\in W_0^{1,2}(B_k)$, we obtain analogs of (\ref{J-11}) and (\ref{J-12}), namely
$$J_k(u)\geq \f{1}{4}\|u\|_{W_0^{1,2}(B_k)}^2-\f{1}{\mu_0}\|f\|_{L^1(V)}^2$$
and
$$\Lambda_k=\inf_{u\in W_0^{1,2}(B_k)}J_k(u)\geq -\f{1}{\mu_0}\|f\|_{L^1(V)}^2.$$
The remaining part of the proof is completely analogous to that of the case $(i)$, and is omitted.\\

{\it Case $(iii)$. $w_{xy}\geq w_0>0$ for all $y$ adjacent to $x$, $\rho\in L^p(V)$ and $f\in L^{{p}/{(p-1)}}(V)$
for some $p\in[1,\infty]$, in particular $f\in L^\infty(V)$ if $p=1$.}

It follows from the Sobolev embedding (Theorem \ref{Soblev-2}) that there exists some constant $C$
depending only on $w_0$, $\mu(O)$, $\|\rho\|_{L^p(V)}$ and $p$ satisfying
$$\|u\|_{L^p(B_k)}\leq C\|u\|_{W_0^{1,2}(B_k)},\quad\forall u\in W_0^{1,2}(B_k).$$
Similar to (\ref{f-bd}), we have
\bna\int_{B_k}fud\mu&\leq& \|u\|_{L^p(B_k)}\|f\|_{L^{\f{p}{p-1}}(V)}\\
&\leq& C\|u\|_{W_0^{1,2}(B_k)}\|f\|_{L^{\f{p}{p-1}}(V)}\\
&\leq&\f{1}{4}\|u\|_{W_0^{1,2}(B_k)}^2+C^2\|f\|_{L^{\f{p}{p-1}}(V)}^2.\ena
As a consequence,  we obtain analogs of (\ref{J-11}) and (\ref{J-12}) as follows:
$$J_k(u)\geq \f{1}{4}\|u\|_{W_0^{1,2}(B_k)}^2-C^2\|f\|_{L^{\f{p}{p-1}}(V)}^2,$$
and
$$\Lambda_k=\inf_{u\in W_0^{1,2}(B_k)}J_k(u)\geq -C^2\|f\|_{L^{\f{p}{p-1}}(V)}^2.$$
Again the remaining part of the proof in this case is completely analogous to that of Case $(i)$,
and thus is omitted.
$\hfill\Box$

\section{Mean field equation}\label{sec5}

In this section, we consider mean field equations. Precisely we prove Theorems \ref{thm4} and \ref{thm5}
by variational method from local to global.

\subsection{The case $g\leq f<0$}

{\it Proof of Theorem \ref{thm4}.} Fix some point $O\in V$. For any $x\in V$, $\rho(x)=\rho(x,O)$ denotes
the distance between $x$ and $O$. For any positive integer $k$, we let
$B_k=\{x\in V:\rho(x)<k\}$, and define a functional $J_k:W_0^{1,2}(B_k)\ra\mathbb{R}$ by
$$J_k(u)=\f{1}{2}\int_{B_k}|\nabla u|^2d\mu+\int_{B_k}fud\mu-\int_{B_k}ge^ud\mu.$$

{\bf Step 1.} {\it For any positive integer $k$, $J_k$ has a lower bound on $W_0^{1,2}(B_k)$.}

Since $g\leq f<0$ and $g\in L^1(V)$, we have also $f\in L^1(V)$. An elementary inequality $e^t\geq 1+t$ for all $t\in \mathbb{R}$
implies that for all $u\in W_0^{1,2}(B_k)$,
\bea\nonumber
J_k(u)&\geq& \int_{B_k}fud\mu-\int_{B_k}ge^ud\mu\\\nonumber
&\geq& \int_{B_k}f(u-e^u)d\mu\\\nonumber
&\geq&\int_{B_k}(-f)d\mu\\\label{J-31}
&=&\int_V(-f)d\mu+o_k(1),
\eea
where $o_k(1)\ra 0$ as $k\ra\infty$. Denoting $c_k=\int_{B_k}(-f)d\mu$, we obtain $J_k(u)\geq c_k$ for all $u\in W_0^{1,2}(B_k)$.

{\bf Step 2.} {\it For any positive integer $k$, there exists a function $u_k\in W_0^{1,2}(B_k)$ such that
\be\label{Juk}J_k(u_k)=\Lambda_k=\inf_{u\in W_0^{1,2}(B_k)}J_k(u).\ee
Moreover $u_k$ satisfies the Euler-Lagrange equation
\be\label{E-L}\le\{\begin{array}{lll}
\Delta u_k=f-ge^{u_k} &{\rm in}& B_k\\[1.5ex]
u_k=0&{\rm on}&\p B_k.
\end{array}\ri.\ee}

Obviously there holds
$$\Lambda_k\leq J_k(0)=\int_{B_k}(-g)d\mu\leq \int_V(-g)d\mu.$$
This together with (\ref{J-31}) gives
\be\label{Lambda-bbd}\|f\|_{L^1(V)}+o_k(1)\leq \Lambda_k\leq \|g\|_{L^1(V)}.\ee
Take a minimizing sequence $(\widetilde{u}_j)\subset W_0^{1,2}(B_k)$ satisfying
\be\label{inf}J_k(\widetilde{u}_j)\ra \Lambda_k=\inf_{u\in W_0^{1,2}(B_k)}J_k(u)\quad{\rm as}\quad j\ra\infty.\ee
For any function $v:V\ra\mathbb{R}$, we write
$$v^+(x)=\le\{\begin{array}{lll}
v(x)&{\rm if}&v(x)> 0\\[1.2ex]
0&{\rm if}& v(x)\leq 0;
\end{array}\ri.\quad v^-(x)=\le\{\begin{array}{lll}
v(x)&{\rm if}&v(x)< 0\\[1.2ex]
0&{\rm if}& v(x)\geq 0.
\end{array}\ri.$$
To see a lower bound of $J_k(\widetilde{u}_j)$, we calculate
\bea\nonumber
J_k(\widetilde{u}_j)&=&\f{1}{2}\int_{B_k}|\nabla \widetilde{u}_j|^2d\mu+\int_{B_k}(f\widetilde{u}_j^+-ge^{\widetilde{u}_j^+})d\mu-\int_{B_k}ge^{\widetilde{u}_j^-}d\mu+\int_{B_k}
f\widetilde{u}_j^-d\mu+\int_{B_k}gd\mu\\
&\geq& \f{1}{2}\int_{B_k}|\nabla \widetilde{u}_j|^2d\mu+\int_{B_k}(f\widetilde{u}_j^+-ge^{\widetilde{u}_j^+})d\mu+\int_{B_k}
f\widetilde{u}_j^-d\mu+\int_{B_k}gd\mu.\label{geq}
\eea
Combining (\ref{inf}) and (\ref{geq}), and noting that $g(x)\leq f(x)<0$, $\widetilde{u}_j^-(x)\leq 0$ for all $x\in B_k$,  we conclude that $(\widetilde{u}_j^-)$ is bounded in $B_k$ with respect to $j$, or equivalently there exists some constant $C$
depending on $k$ such that
\be\label{u-}|\widetilde{u}_j^-(x)|\leq C\quad{\rm for\,\,all}\quad x\in B_k.\ee
Note also that
$$f\widetilde{u}_j^+-ge^{\widetilde{u}_j^+}\geq f\widetilde{u}_j^+-fe^{\widetilde{u}_j^+}\geq -\f{f}{2}{(\widetilde{u}_j^+)}^2, $$
which together with (\ref{inf}) and (\ref{geq}) leads to
\be\label{u+}\widetilde{u}_j^+(x)\leq C\quad{\rm for\,\,all}\quad x\in B_k,\ee
where $C$ is some constant depending on $k$. It follows from (\ref{u-}) and (\ref{u+}) that
$(u_j)$ is uniformly bounded in $B_k$ with respect to $j$. Hence there exist a subsequence of $(\widetilde{u}_j)$, which is still denoted by $(\widetilde{u}_j)$,
and a function ${u}_k\in W_0^{1,2}(B_k)$ such that
$\widetilde{u}_j$ converges to ${u_k}$ uniformly in $B_k$ as $j\ra\infty$. This together with (\ref{inf}) immediately leads to
(\ref{Juk}).
By a straightforward calculation, ${u}_k$ satisfies the Euler-Lagrange equation (\ref{E-L}).\\

{\bf Step 3.} {\it For any finite set $A\subset V$, $(u_k)$ is uniformly bounded in $A$.}

Let $A$ be a finite subset of $V$. An obvious analog of (\ref{geq}) reads
\bna\nonumber
J_k(u_k)&\geq&\f{1}{2}\int_{B_k}|\nabla u_k|^2d\mu+\int_{B_k}(fu_k^+-ge^{u_k^+})d\mu+\int_{B_k}fu_k^-d\mu+\int_{B_k}gd\mu\\
&\geq&\int_{A}(fu_k^+-ge^{u_k^+})d\mu+\int_{A}fu_k^-d\mu+\int_{B_k}gd\mu,
\ena
provided that $k$ is sufficiently large. As a consequence, one derives
\be\label{upb-1}\max_{x\in A}|u_k^-(x)|\leq \f{J_k(u_k)-\int_{B_k}gd\mu}{\min_{x\in A}\mu(x)|f(x)|};\quad
\max_{x\in A}u_k^+(x)\leq \sqrt{\f{2J_k(u_k)-2\int_{B_k}gd\mu}{\min_{x\in A}\mu(x)|f(x)|}}.\ee
Combining (\ref{Juk}),  (\ref{Lambda-bbd}) and (\ref{upb-1}), we conclude that there exists some constant $C$ depending only on
$h$, $g$, $\mu$ and $A$
such that
$$\max_{x\in A}|u_k(x)|\leq C.$$

{\bf Step 4.} {\it There exists a subsequence of $(u_k)$, which is still denoted by $(u_k)$, and a function
$u^\ast:V\ra\mathbb{R}$ such that $(u_k)$ converges to $u^\ast$ locally uniformly in $V$. Moreover, $u^\ast$ is a solution of
the equation (\ref{KZ-equ}).}

 By Step 3, $(u_k)$ is uniformly bounded in $B_1$. Hence there exists a subsequence of $(u_k)$, which is written as
 $(u_{1,k})$, and a function $u_1^\ast$ such that $u_{1,k}$ converges to $u_1^\ast$ in $B_1$.
 By Step 3 again, $(u_{1,k})$ is uniformly bounded in $B_2$. Then there would exist a subsequence of $(u_{1,k})$, which is written as
 $(u_{2,k})$, and a function $u_2^\ast$ such that $u_{2,k}$ convergence to $u_2^\ast$ uniformly in $B_2$. Obviously $u_{2}^\ast=u_1^\ast$
 on $B_1$. Repeating this process, one finds a diagonal subsequence $(u_{k,k})$, which is still denoted by
 $(u_{k})$, and a function $u^\ast:V\ra\mathbb{R}$ such that for any finite set $A\subset V$, $(u_k)$ converges to $u^\ast$ uniformly in $A$.
 For any fixed $x\in V$, passing to the limit $k\ra\infty$ in (\ref{E-L}), we obtain
 $$\Delta u^\ast(x)=f(x)-g(x)e^{u^\ast(x)}.$$
 This ends the final step and completes the proof of the theorem. $\hfill\Box$

 \subsection{The case $g>0$}

 {\it Proof of Theorem \ref{thm5}.} Fix some point $O\in V$. For any $x\in V$, $\rho(x)=\rho(x,O)$ denotes the distance between
 $x$ and $O$. Let $B_k=\{x\in V:\rho(x)<k\}$, $W_0^{1,2}(B_k)$ be the Sobolev space including all
 functions $u$ satisfying $u=0$ on $\p B_k$, with the norm
 $$\|u\|_{W_0^{1,2}(B_k)}=\le(\int_{B_k}(|\nabla u|^2+hu^2)d\mu\ri)^{1/2},$$
 where $h(x)\geq a_0>0$, $\mu(x)\geq \mu_0>0$ for all $x\in V$.
  Define a functional $J_k:W_0^{1,2}(B_k)\ra\mathbb{R}$ by
 \be\label{funct}J_k(u)=\f{1}{2}\int_{B_k}(|\nabla u|^2+hu^2)d\mu+\int_{B_k}fud\mu-\log\int_{B_k}ge^{u}d\mu.\ee
  Since $f\in L^q(V)$ for some $q$ with $1\leq q\leq 2$, we have
 by the Sobolev embedding (Theorem \ref{Soblev-1}),
 \be\label{h-b}\le|\int_{B_k}fud\mu\ri|\leq \|f\|_{L^q(V)}\|u\|_{L^p(B_k)}\leq C\|f\|_{L^q(V)}\|u\|_{W_0^{1,2}(B_k)}\ee
 for some constant $C$ depending only on $\mu_0$, $a_0$ and $q$, where $1/p+1/q=1$. Since
 $$\|v\|_{L^\infty(B_k)}\leq \f{1}{\sqrt{\mu_0a_0}}\|v\|_{W_0^{1,2}(B_k)},\quad\forall v\in W_0^{1,2}(B_k),$$
 there holds for any $\epsilon>0$,
 $$e^{u}\leq e^{\f{u^2}{4\epsilon \|u\|_{W_0^{1,2}(B_k)}^2}+\epsilon\|u\|_{W_0^{1,2}(B_k)}^2}\leq e^{\f{1}{4\epsilon\mu_0a_0}+\epsilon
 \|u\|_{W_0^{1,2}(B_k)}^2}.$$
 It then follows that
 \be\label{log}\log\int_{B_k}ge^ud\mu\leq \log\|g\|_{L^1(V)}+\f{1}{4\epsilon\mu_0a_0}+\epsilon
 \|u\|_{W_0^{1,2}(B_k)}^2.\ee
 Note that $\|g\|_{L^1(V)}>0$, since $g\geq 0$ but $g\not\equiv 0$. Inserting (\ref{h-b}) and (\ref{log}) into (\ref{funct}), we obtain
 $$J_k(u)\geq \le(\f{1}{2}-\epsilon\ri)\|u\|_{W_0^{1,2}(B_k)}^2-\f{1}{4}\|u\|_{W_0^{1,2}(B_k)}^2
 -C^2\|f\|_{L^q(V)}^2-\log\|g\|_{L^1(V)}-\f{1}{4\epsilon\mu_0a_0}.$$
 Choosing $\epsilon=1/8$, we immediately have for any $u\in W_0^{1,2}(B_k)$,
 \be\label{J-low}J_k(u)\geq \f{1}{8}\|u\|_{W_0^{1,2}(B_k)}^2-C^2\|f\|_{L^q(V)}^2-\log\|g\|_{L^1(V)}-\f{2}{\mu_0a_0}.\ee
 Hence $J_k$ has a lower bound in $W_0^{1,2}(B_k)$. Take a minimizing sequence $(\widetilde{u}_j)\subset W_0^{1,2}(B_k)$ such that
 \be\label{seq}J_k(\widetilde{u}_j)\ra \Lambda_k=\inf_{u\in W_0^{1,2}(B_k)}J_k(u)\quad{\rm as}\quad j\ra\infty.\ee
 Since $g\geq 0$ and there exists some $x_0\in V$ such that $g(x_0)>0$, there holds
 $$\mu(x_0)g(x_0)\leq \int_{B_k}gd\mu,$$
 and thus
 \be\label{A-bd}\Lambda_k\leq J_k(0)=-\log\int_{B_k}gd\mu\leq -\log(\mu(x_0)g(x_0)).\ee
 Combining (\ref{J-low}), (\ref{seq}) and (\ref{A-bd}), we have
 $$\|\widetilde{u}_j\|_{W_0^{1,2}(B_k)}\leq C$$
 for some constant $C$ independent of $k$. Hence there exists a subsequence of $(\widetilde{u}_j)$, which is still
 denoted by $(\widetilde{u}_j)$,
 and a function $u_k\in W_0^{1,2}(B_k)$ such that $(\widetilde{u}_j)$ converges to $u_k$ uniformly in $B_k$ as $j\ra\infty$. It is easy to see that
 $u_k$ is a minimizer of $J_k$, or equivalently
 $$J_k(u_k)=\Lambda_k=\inf_{u\in W_0^{1,2}(B_k)}J_k(u).$$
 Moreover $u_k$ satisfies the Euler-Lagrange equation
 \be\label{EL-2}\le\{
 \begin{array}{lll}
 -\Delta u_k+hu_k=\f{1}{\gamma_k}ge^{u_k}-f\,\,\,{\rm in}\,\,\, B_k\\[1.5ex]
 u_k\in W_0^{1,2}(B_k),\,\,\,\gamma_k=\int_{B_k}ge^{u_k}d\mu.
 \end{array}
 \ri.\ee
  Since $(\Lambda_k)$ is bounded due to (\ref{J-low}) and (\ref{A-bd}), we conclude that
  \be\label{W-b}\|u_k\|_{W_0^{1,2}(B_k)}\leq C\ee
  for some constant $C$ independent of $k$. Using the same argument as Step 4 of the proof
  of Theorem \ref{thm4}, one easily extracts a subsequence of
  $u_k$, which is still denoted by $u_k$, and finds some function $u^\ast$ such that
  $(u_k)$ converges to $u^\ast$ locally uniformly in $V$. In view of (\ref{W-b}), the Sobolev embedding theorem
  (Theorem \ref{Soblev-1}) implies
  \be\label{uk-up}\|u_k\|_{L^\infty(B_k)}\leq \f{1}{\sqrt{\mu_0}}\|u_k\|_{W_0^{1,2}(B_k)}\leq C.\ee
  This immediately leads to
  $$e^{-C}\|g\|_{L^1(B_k)}\leq \gamma_k\leq e^{C}\|g\|_{L^1(B_k)},$$
  where $\gamma_k$ is given as in (\ref{EL-2}).
  Then up to a subsequence, $\gamma_k$ converges to some number $\gamma^\ast$ with
  \be\label{ga-m}e^{-C}\|g\|_{L^1(V)}\leq \gamma^\ast\leq e^{C}\|g\|_{L^1(V)}.\ee
  It follows from (\ref{EL-2}) and (\ref{ga-m}) that
  \be\label{equ-1}-\Delta u^\ast+hu^\ast=\f{1}{\gamma^\ast}ge^{u^\ast}-f\quad{\rm in}\quad V.\ee

  We now prove  \be\label{ga-repre}\gamma^\ast=\int_Vge^{u^\ast}d\mu.\ee
  On one hand, for any fixed $\ell>1$, there holds
  $$\int_{B_\ell}ge^{u^\ast}d\mu=\lim_{k\ra\infty}\int_{B_\ell}ge^{u_k}d\mu\leq
  \lim_{k\ra\infty}\int_{B_k}ge^{u_k}d\mu=\gamma^\ast,$$
  which leads to
  \be\label{inequ-1}\int_Vge^{u^\ast}d\mu\leq \gamma^\ast.\ee
  On the other hand, in view of (\ref{uk-up}) and the assumption $g\in L^1(V)$, for any $\eta>0$, there would exist a sufficiently large $\ell_0>1$ such that
  if $\ell\geq \ell_0$, then
  \be\label{lim-1}\int_{B_k}ge^{u_k}d\mu\leq \eta+\int_{B_\ell}ge^{u_k}d\mu.\ee
  Indeed, (\ref{uk-up}) and $g\in L^1(V)$ lead to
  $$\int_{B_k\setminus B_\ell}ge^{u_k}d\mu\leq e^C\int_{V\setminus B_\ell}gd\mu=o_\ell(1),$$
  where $o_\ell(1)\ra 0$ as $\ell\ra\infty$.
  Thus (\ref{lim-1}) is satisfied. Passing to the limit $k\ra\infty$ first, and then $\ell\ra\infty$ in (\ref{lim-1}), we obtain
  $$\gamma^\ast\leq \eta+\int_Vge^{u^\ast}d\mu.$$
  Since $\eta>0$ is arbitrary, there must hold
  \be\label{l-3}\gamma^\ast\leq \int_Vge^{u^\ast}d\mu.\ee
  Hence (\ref{ga-repre}) follows from (\ref{inequ-1}) and (\ref{l-3}) immediately.
  Combining (\ref{ga-repre}) and (\ref{equ-1}), we conclude that $u^\ast$ is a solution of
  $$\le\{
 \begin{array}{lll}
 -\Delta u^\ast+hu^\ast=\f{1}{\gamma^\ast}ge^{u^\ast}-f\,\,\,{\rm in}\,\,\, V\\[1.5ex]
 \gamma^\ast=\int_{V}ge^{u^\ast}d\mu.
 \end{array}
 \ri.$$

 Since $u_k$ is naturally viewed as a function on $V$, using the same argument as the proof of (\ref{uast}), we
 conclude  from (\ref{W-b}) and (\ref{H-norm}) that $u^\ast\in\mathscr{H}$.
   This completes the proof of the theorem. $\hfill\Box$

 \section{Yamabe equation}\label{sec6}

 In this section, using the mountain-pass theorem due to Ambrosetti-Rabinowitz \cite{Ambrosetti-Rabinowitz},
 we prove the existence of nontrivial solutions to the Yamabe equation (\ref{Yamabe}).
 The key estimate is the Sobolev embedding theorem. In \cite{GLY3}, we have used Theorem \ref{Soblev-1} under the assumption
 (\ref{mu0}). Here we shall apply Theorem \ref{Soblev-2} to the mountain-pass theory.
 Our assumptions on the locally finite graph are $w_{xy}\geq w_0>0$ for all $y$ adjacent to $x$, and
 $$\int_V\rho^pd\mu=\sum_{x\in V}\mu(x)\rho^p(x)<+\infty$$
 for some $p>2$,
 where $\rho(x)=\rho(x,O)$ denotes the distance between $x$ and $O$. It seems that Theorem \ref{Soblev-2} has a lot of room for
 improvement.\\

 To begin with, we have the following compactness embedding
 for $\mathscr{H}$, where $\mathscr{H}$ is a Hilbert space defined as in (\ref{H}).
 \begin{lemma}\label{compactembedding}
 If $h\geq a_0>0$ and $1/h\in L^1(V)$, then $\mathscr{H}$ is embedded in $L^q(V)$ compactly for all $1\leq q<p$;
 If $h\geq a_0>0$ and $h(x)\ra+\infty$ as $\rho(x)\ra+\infty$, then $\mathscr{H}$ is embedded in $L^q(V)$
 compactly for all $2\leq q<p$.
 \end{lemma}

 \proof Suppose $(u_k)$ is a bounded sequence in $\mathscr{H}$, namely
 \be\label{H-bd}\|u_k\|_{\mathscr{H}}^2=\int_V(|\nabla u_k|^2+hu_k^2)d\mu\leq C.\ee
 Since the Hilbert space $\mathscr{H}$ is reflexive, there exists some function $u\in \mathscr{H}$ such that up to
 a subsequence,
 $(u_k)$ converges to $u$ weakly in $\mathscr{H}$, locally uniformly in $V$. If $1/h\in L^1(V)$, then for any
 $\epsilon>0$, there exists some $\ell>1$ such that
 $$\int_{V\setminus B_\ell}\f{1}{h}d\mu<\epsilon^2.$$
 Moreover, there holds
 \bna
 \int_V|u_k-u|d\mu&\leq&\int_{B_\ell}|u_k-u|d\mu+\le(\int_{V\setminus B_\ell}\f{1}{h}d\mu\ri)^{1/2}
 \le(\int_{V\setminus B_\ell}h|u_k-u|^2d\mu\ri)^{1/2}\\
 &\leq& C\epsilon+o_k(1).
 \ena
 This immediately implies
 \be\label{L1-conv}\lim_{k\ra\infty}\|u_k-u\|_{L^1(V)}=0.\ee
 For any $q\in(1,p)$, there exists a unique $\lambda\in (0,1)$ such that $q=\lambda+(1-\lambda)p$. By the H\"older inequality, (\ref{H-bd})
 and Theorem \ref{Soblev-2},
 \bna
 \int_V|u_k-u|^qd\mu&\leq& \le(\int_V|u_k-u|d\mu\ri)^\lambda\le(\int_V|u_k-u|^pd\mu\ri)^{1-\lambda}\\
 &\leq & C\le(\int_V|u_k-u|d\mu\ri)^\lambda,
 \ena
 which together with (\ref{L1-conv}) leads to
 \be\label{Lq}\lim_{k\ra\infty}\|u_k-u\|_{L^q(V)}=0.\ee

 If $h(x)\ra+\infty$ as $\rho(x)\ra\infty$, then for any $\epsilon>0$, there exists some $\ell_1>1$ such that
 $$\label{L2}h(x)\geq \f{C}{\epsilon}\quad{\rm for\,\,all}\quad x\in V\setminus B_{\ell_1}.$$
 As a consequence
 \bna
 \int_V|u_k-u|^2d\mu&=&\int_{B_{\ell_1}}|u_k-u|^2d\mu+\int_{V\setminus B_{\ell_1}}|u_k-u|^2d\mu\\
 &\leq&\f{\epsilon}{C}\int_{V\setminus B_{\ell_1}}h|u_k-u|^2d\mu+o_k(1).
 \ena
 This implies that
 \be\label{62}\lim_{k\ra\infty}\|u_k-u\|_{L^2(V)}=0.\ee
 Using the same argument as in the proof of (\ref{Lq}), we obtain from (\ref{62}) that
 $$\lim_{k\ra\infty}\|u_k-u\|_{L^q(V)}=0\quad{\rm for\,\,all}\quad 2< q<p.$$
 This ends the proof of the lemma.$\hfill\Box$\\

 Let $f$ be a function of one variable defined by
 \be\label{f}f(s)=|s|^{q-2}s,\quad s\in \mathbb{R}\ee
 and $F$ be its primitive function, namely
 \be\label{F}F(s)=\int_0^sf(t)dt=\f{1}{q}|s|^q,\quad s\in \mathbb{R}.\ee
 Obviously $sf(s)=qF(s)$ for all $s\in\mathbb{R}$.
 Define a functional $J:\mathscr{H}\ra\mathbb{R}$ by
 \be\label{funct-yamabe}J(u)=\f{1}{2}\int_V(|\nabla u|^2+hu^2)d\mu-\int_{V}F(u)d\mu.\ee
 \begin{lemma}\label{PS}
 Assume $q\in (2,p)$, $f$, $F$ and $J$ are defined as in (\ref{f}), (\ref{F}) and (\ref{funct-yamabe}) respectively. Then
 for any $c\in\mathbb{R}$, $J$ satisfies the $(PS)_c$ condition. Precisely, if for any
 sequence $(u_k)\subset \mathscr{H}$ with $J(u_k)\ra c$ and $J^\prime(u_k)\ra 0$, then up to a subsequence,
 $(u_k)\ra u$ in $\mathscr{H}$ for some function $u\in \mathscr{H}$.
 \end{lemma}
 \proof Since $(u_k)\subset \mathscr{H}$, $J(u_k)\ra c$ and $J^\prime(u_k)\ra 0$, we have
 \bea\label{J-c}
 &&\f{1}{2}\|u_k\|_{\mathscr{H}}^2-\int_VF(u_k)d\mu=c+o_k(1)\\\label{J'c}
 &&\langle u_k,\phi\rangle_{\mathscr{H}}-\int_Vf(u_k)\phi d\mu=o_k(1)\|\phi\|_{\mathscr{H}},\quad\forall \phi\in\mathscr{H}.
 \eea
 Taking $\phi=u_k$ in (\ref{J'c}) and noting that $u_k(x)f(u_k(x))=qF(u_k(x))$ for all $x\in V$, we obtain
 \be\label{comb}\f{q}{2}\|u_k\|_{\mathscr{H}}^2-qc+o_k(1)=\|u_k\|_{\mathscr{H}}^2+o_k(1)\|u_k\|_{\mathscr{H}}.\ee
 Since $2<q<p$, (\ref{comb}) implies that $(u_k)$ is bounded in $\mathscr{H}$. By Lemma \ref{compactembedding}, there
 exist a subsequence of $(u_k)$, which is still denoted by $(u_k)$, and some function $u\in \mathscr{H}$ such that
 \be\label{qconv}\lim_{k\ra\infty}\int_V|u_k-u|^qd\mu=0.\ee
 One calculates
 \bea\nonumber
 \int_V|F(u_k)-F(u)|d\mu&=&\int_V|f(\xi_k)||u_k-u|d\mu\\\nonumber
 &\leq&\int_V(|u_k|^{q-1}+|u|^{q-1})|u_k-u|d\mu\\\label{F-con}
 &\leq&C\|u_k-u\|_{L^q(V)},
 \eea
 where Theorem \ref{Soblev-2} is used, $C$ is a constant independent of $k$, and $\xi_k$ lies between $u_k$ and $u$. Combining (\ref{qconv}) and (\ref{F-con}), we obtain
 \be\label{int-Fk}\lim_{k\ra\infty}\int_VF(u_k)d\mu=\int_VF(u)d\mu.\ee
 In the same way,
 \bna
 \le|\int_Vf(u_k)(u_k-u)d\mu\ri|&\leq&\le(\int_V|f(u_k)|^{\f{q}{q-1}}d\mu\ri)^{1-1/q}\le(\int_V|u_k-u|^qd\mu\ri)^{1/q}\\
 &\leq&\|u_k\|_{L^q(V)}^{q-1}\|u_k-u\|_{L^q(V)}\\
 &\leq&C\|u_k-u\|_{L^q(V)}.
 \ena
 As a consequence
 \be\label{fkuk}\lim_{k\ra\infty}\int_Vf(u_k)(u_k-u)d\mu=0.\ee
 Taking $\phi=u_k-u$ in (\ref{J'c}) and noting (\ref{fkuk}), we obtain
 \be\label{formula-1}\langle u_k,u_k-u\rangle_{\mathscr{H}}=o_k(1).\ee
 Since up to a subsequence, $u_k\rightharpoonup u$ weakly in $\mathscr{H}$, it follows that
 \be\label{formula-2}\langle u,u_k-u\rangle_{\mathscr{H}}=o_k(1).\ee
 Combining (\ref{formula-1}) and (\ref{formula-2}), we conclude that $(u_k)$ converges to
 $u$ in $\mathscr{H}$. In view of (\ref{J-c}), (\ref{J'c}), (\ref{int-Fk}) and (\ref{fkuk}), we have
 $$J(u)=c,\quad J^\prime(u)=0.$$
 This ends the proof of the lemma. $\hfill\Box$\\

 {\it Proof of Theorem \ref{thm6}.} Let $J\in C^1(\mathscr{H}, \mathbb{R})$ be the functional defined as in (\ref{funct-yamabe}). We claim that $J$
 satisfies $(H_1)$ $J(0)=0$; $(H_2)$ for some $\delta>0$, $\inf_{\|u\|_{\mathscr{H}}=\delta}J(u)>0$; $(H_3)$ $J(v)<0$ for some $v\in \mathscr{H}$
 with $\|v\|_{\mathscr{H}}>\delta$. Firstly, $(H_1)$ is obvious. Secondly, to see $(H_2)$, we have by Lemma \ref{compactembedding},
 \bna\nonumber
 J(u)&\geq&\f{1}{2}\|u\|_{\mathscr{H}}^2-\f{1}{q}\int_V|u|^qd\mu\\\label{SL}
 &\geq&\f{1}{2}\|u\|_{\mathscr{H}}^2-C\|u\|_{\mathscr{H}}^q
 \ena
 for some constant $C$ depending on $q$. Hence, if $\|u\|_{\mathscr{H}}=\delta$ for sufficiently small $\delta>0$, there holds
 $J(u)\geq C>0$ for some constant $C$ depending on $q$ and $\delta$. This confirms $(H_2)$. Finally, to see $(H_3)$, we take a function
 $$u_0(x)=\le\{\begin{array}{lll}
 1,&&x=O\\[1.5ex]
 0,&&x\not=O
 \end{array}\ri.$$
 for some fixed point $O\in V$. It then follows that
 \bna
 J(tu_0)&=&\f{t^2}{2}\|u_0\|_{\mathscr{H}}^2-\f{t^q}{q}\int_Vu_0^qd\mu\\
 &\ra& -\infty\quad{\rm as}\quad t\ra+\infty,
 \ena
 since $2<q<p$.
 If we choose $v=tu_0$ for sufficiently large $t>0$, then $J(v)<0$ and $(H_3)$ holds.

 Let $$\label{min-max}c=\min_{\gamma\in\Gamma}\max_{u\in\gamma}J(u),$$
 where $\Gamma=\{\gamma|\gamma:[0,1]\ra\mathscr{H}\,\,{\rm is\,\,a\,\,} C^1\,\,{\rm curve\,\, with}\,\,\gamma(0)=0,\gamma(1)=v \}$. Clearly
 $\label{c}0<c<+\infty$. In view of Lemma \ref{PS}, applying the mountain-pass theorem due to Ambrosetti-Rabinowitz \cite{Ambrosetti-Rabinowitz}, we conclude that $c$ is a critical value of $J$. In particular, there exists some $u\in \mathscr{H}$ such that
 $J(u)=c$, $J^\prime(u)=0$. Clearly $u\not\equiv 0$, and $u$ satisfies the Euler-Lagrange equation
 (\ref{Yamabe}).
 $\hfill\Box$\\

 {\bf Acknowledgements.} We thank the reviewers for their careful reading and valuable comments. Yong Lin is supported by the NSFC (Grant No. 12071245).
 Yunyan Yang is supported by the  NSFC (Grant No. 11721101) and
     the National Key Research and Development Project SQ2020YFA070080.
 Both of the two authors are supported by the NSFC (Grant No. 11761131002).

\end{document}